# Some vaccination strategies for the SEIR epidemic model. Preliminary results

M. De la Sen*, A. Ibeas ** , S. Alonso-Quesada*

Department of Electricity and Electronics. Faculty of Science and Technology

University of the Basque Country. PO. Box 644- Bilbao. Spain

** Department of Telecommunications and Systems Engineering

Universitat Autònoma de Barcelona, 08193- Bellaterra, Barcelona, Spain

**Abstract**. This paper presents a vaccination-based control strategy for a SEIR (susceptible plus infected plus infectious plus removed populations) propagation disease model. The model takes into account the total population amounts as a refrain for the illness transmission since its increase makes more difficult contacts among susceptible and infected. The control objective is the asymptotically tracking of the removed-by-immunity population to the total population while achieving simultaneously the remaining population (i.e. susceptible plus infected plus infectious) to asymptotically tend to zero.

## 1. Introduction

Important control problems nowadays related to Life Sciences are the control of ecological models are, for instance, those of population evolution ( Beverton-Holt model, Hassell model, Ricker model etc.) via the online adjustment of the species environment carrying capacity, that of the population growth or that of the regulated harvesting quota as well as the disease propagation via vaccination control. In a set of papers, several variants and generalizations of the Beverton-Holt model (standard time–invariant, time-varying parameterized, generalized model or modified generalized model) have been investigated at the levels of stability, cycle- oscillatory behavior, permanence and control through the manipulation of the carrying capacity (see, for instance, [1-5]). The design of related control actions has been proved to be important in those papers at the levels, for instance, of aquaculture exploitation or plague fighting. On the other hand, the literature about epidemic mathematical models is exhaustive in many books and papers. A non-exhaustive list of references is given in this manuscript, cf. [6-14] (see also the references listed therein). The sets of models are described in [6-7]. Those models have also two major variants, namely, the so-called "pseudo-mass action models", where the total population is not taken into account as a relevant disease contagious factor and the so-called "true-mass action models", where the total population is more realistically considered as an inverse factor of the disease transmission rates. There are many variants of the above models, for instance, including vaccination of different kinds: constant [8], impulsive [12], discrete – time etc., incorporating point or distributed delays [12-13], oscillatory behaviors [14] etc. . Some nonlinear models have been proven to be also useful to deal with kinetic modelling because of their rich dynamics issues (see, for instance, [15]). On the other hand, variants of such models become considerably simpler for the illness transmission among plants [6-7]. In this paper, a continuous-time vaccination control strategy is given for a SEIR epidemic model which takes directly the total, infectious and removed– by-immunity numbers to design the vaccination strategy. It is assumed that the total population remains uniformly bounded through time while being nonnegative as they are all the partial populations of susceptible, infected, infectious and immune. Thus, the disease transmission is not critical, and the SEIR – model is of the above mentioned true-mass action type. Note that although all the partial populations and the total one are all nonnegative for all time in the real problem under study, the property has to be guaranteed for the mathematical SEIR- model (1)-(4) as well.

## 2. A true-mass action SEIR epidemic model

Let S (t) be the "susceptible" population of infection at time t, E (t) the " infected" ( i.e. those which incubate the illness but do not still have any symptoms) at time t, I (t ) is the " infectious" ( or "infective") population at time t, and R (t) is the "removed by immunity" ( or " immune") population at time t. Consider the true-mass type SEIR-type epidemic model:

$$\dot{S}(t) = -\mu S(t) + \omega R(t) - \beta \frac{S(t)I(t)}{N(t)} + \nu N(t)(1 - V(t))  \quad (1)$$

$$\dot{E}(t) = \beta \frac{S(t)I(t)}{N(t)} - (\mu + \sigma)E(t)  \quad (2)$$

$$\dot{I}(t) = -(\mu + \gamma)I(t) + \sigma E(t)  \quad (3)$$

$$\dot{R}(t) = -(\mu + \omega)R(t) + \gamma(1 - \rho)I(t) + \nu N(t)V(t)  \quad (4)$$

subject to initial conditions $S_0 = S(0) \geq 0$, $E_0 = E(0) \geq 0$, $I_0 = I(0) \geq 0$ and $R_0 = R(0) \geq 0$ under the vaccination function $V : \mathbf{R}_{0+} \to \mathbf{R}_+$. The vaccination control is either the vaccination function itself or some appropriate four dimensional vector depending on it defined "ad –hoc" for some obtained equivalent representation of the SEIR- model as a dynamic system. In the above SEIR – model, N is the total population, μ is the rate of deaths from causes unrelated to the infection, ω is the rate of losing immunity, β is the transmission constant (with the total number of infections per unity of time at time t being $\beta \frac{S(t)I(t)}{N(t)}$), $\sigma^{-1}$ and $\gamma^{-1}$ are, respectively, the average durations of the latent and infective periods. All the above parameters are nonnegative. The parameter $\omega$ is the rate of immunity lost since it makes the susceptible to increase and then the immune to decrease. The usual simplified SEIR- model is obtained with the standard model particular case $\nu = \mu$, where ν is related to the vaccination of newborns and the parameter ρ of dead caused by the disease is made zero in the standard version of the model, [7]. In that standard particular case of the model, one gets:

$$\dot{N}(t) = \dot{S}(t) + \dot{E}(t) + \dot{I}(t) + \dot{R}(t) = \mu(N(t) - S(t) - E(t) - I(t) - R(t)) = 0 \; ; \; \forall t \in \mathbf{R}_{0+}$$

$$\Rightarrow N(t) := S(t) + E(t) + I(t) + R(t) = N(0) = N_0 = N > 0$$

If $\nu > \mu$ then the new-born lost of maternal immunity is considered in the model what makes the susceptibility rate to increase and, perhaps, to increase mortality by cause of the disease. If $\nu < \mu$ then there is a rapid vaccination action to newborns to make them to be removed from the susceptible. The parameter $\rho \in (0,1]$ is the per- capita probability of dying from the infection. If either both $\nu \neq \mu$ and $\rho = 0$ or, alternatively, both $\nu = \mu$ and $\rho \neq 0$ then N(t) is not constant through time. If $\nu = \mu$ and $\rho = 0$ then N(t) = N; $\forall t \in \mathbf{R}_{0+}$ is constant through time. If, in addition, the vaccination function is identically zero then one obtains the standard true-mass type SEIR epidemic model with constant population through time as particular case of the model (1)-(5). If $I(t) = \dfrac{(v-\mu)N(t)}{\rho \gamma}$ occurs eventually on a set of zero measure then the total population varies through time as obtained by correspondingly summing- up both sides of (1)-(4) yielding $\dot{N}(t) = (v-\mu)N(t) - \rho \gamma I(t)$; $\forall t \in R_{0+}$. Eqs. (1)-(4) describe the state evolution of the SEIR- model with state vector $x(t) = (S(t), E(t), I(t), R(t))^T$, output $y(t) \equiv N(t)$ and control $V(t)$ as the following nonlinear dynamic system which involves the nonlinear term $\dfrac{S(t)I(t)}{N(t)}$:

$$\dot{x}(t) = A_i(x(t), t)x(t) + b_a y(t)V(t) + \zeta y(t) \; ; \; i=1,2,3,4 \tag{5}$$

$$= A_i(x(t), t)x(t) + v N(t) u(t) \tag{6}$$

$$= A_i^0(x(t), t)x(t) + b_a y(t) V(t) \tag{7}$$

$$y(t) \equiv N(t) = c^T x(t) \tag{8}$$

where

$$b_a = (-v, 0, 0, v)^T \; ; \; \zeta := v e_1 = (v, 0, 0, 0)^T \; ; \; c = (1, 1, 1, 1)^T$$

$$u(t) = (1-V(t), 0, 0, V(t))^T = (1-V(t))e_1 + e_4 V(t) = e_1 + (e_4 - e_1)V(t) \tag{9}$$

with $e_i$ $(1 \le i \le 4)$ being the unit fourth dimensional Euclidean vector with i-th component being one and the matrix of dynamics in (5)-(8) being any of the subsequent ones:

$$A_1(x(t), t) = \begin{bmatrix} -\left(\mu + \dfrac{\beta I(t)}{N(t)}\right) & 0 & 0 & \omega \\ \dfrac{\beta I(t)}{N(t)} & -(\mu+\sigma) & 0 & 0 \\ 0 & \sigma & -(\mu+\gamma) & 0 \\ 0 & 0 & \gamma(1-\rho) & -(\mu+\omega) \end{bmatrix}$$

$$A_2(x(t), t) = \begin{bmatrix} -\mu & 0 & -\dfrac{\beta S(t)}{N(t)} & \omega \\ 0 & -(\mu+\sigma) & \dfrac{\beta S(t)}{N(t)} & 0 \\ 0 & \sigma & -(\mu+\gamma) & 0 \\ 0 & 0 & \gamma(1-\rho) & -(\mu+\omega) \end{bmatrix}$$

$$A_3(x(t),t) = \begin{bmatrix} -\mu & 0 & -\frac{\beta S(t)}{N(t)} & \omega \\ \frac{\beta I(t)}{N(t)} & -(\mu+\sigma) & 0 & 0 \\ 0 & \sigma & -(\mu+\gamma) & 0 \\ 0 & 0 & \gamma(1-\rho) & -(\mu+\omega) \end{bmatrix}$$

$$A_4(x(t),t) = \begin{bmatrix} -\left(\mu+\frac{\beta I(t)}{N(t)}\right) & 0 & 0 & \omega \\ 0 & -(\mu+\sigma) & \frac{\beta S(t)}{N(t)} & 0 \\ 0 & \sigma & -(\mu+\gamma) & 0 \\ 0 & 0 & \gamma(1-\rho) & -(\mu+\omega) \end{bmatrix}$$

$$A_1^0(x(t),t) = A_1(x(t),t) + \zeta c^T c^T x(t) = \begin{bmatrix} -\left(\mu+\frac{\beta I(t)}{N(t)}-v\right) & v & v & \omega+v \\ \frac{\beta I(t)}{N(t)} & -(\mu+\sigma) & 0 & 0 \\ 0 & \sigma & -(\mu+\gamma) & 0 \\ 0 & 0 & \gamma(1-\rho) & -(\mu+\omega) \end{bmatrix}$$

$$A_2^0(x(t),t) = A_2(x(t),t) + \zeta c^T c^T x(t) = \begin{bmatrix} v-\mu & v & v-\frac{\beta S(t)}{N(t)} & \omega+v \\ 0 & -(\mu+\sigma) & \frac{\beta S(t)}{N(t)} & 0 \\ 0 & \sigma & -(\mu+\gamma) & 0 \\ 0 & 0 & \gamma(1-\rho) & -(\mu+\omega) \end{bmatrix}$$

$$A_3^0(x(t),t) = A_3(x(t),t) + \zeta c^T c^T x(t) = \begin{bmatrix} v-\mu & v & v-\frac{\beta S(t)}{N(t)} & \omega+v \\ \frac{\beta I(t)}{N(t)} & -(\mu+\sigma) & 0 & 0 \\ 0 & \sigma & -(\mu+\gamma) & 0 \\ 0 & 0 & \gamma(1-\rho) & -(\mu+\omega) \end{bmatrix}$$

$$A_4^0(x(t),t) = A_4(x(t),t) + \zeta c^T c^T x(t) = \begin{bmatrix} -\left(\mu+\frac{\beta I(t)}{N(t)}-v\right) & v & v & \omega+v \\ 0 & -(\mu+\sigma) & \frac{\beta S(t)}{N(t)} & 0 \\ 0 & \sigma & -(\mu+\gamma) & 0 \\ 0 & 0 & \gamma(1-\rho) & -(\mu+\omega) \end{bmatrix}$$

(10)

### 3. Positivity

Note that a necessary condition for a forced linear time-invariant system to be positive in the sense that its state components and output are nonnegative for all time is that its matrix of dynamics be a Metzler one, namely, all its off-diagonal entries are nonnegative. This is equivalent to say that the corresponding unforced system (i.e. that obtained under identically zero control) has to be positive as a necessary condition for the positivity of the forced one. Note that the fundamental matrix whose infinitesimal generator is a Metzler matrix is a positive $C_0$-semigroup, namely, all its entries are nonnegative. Now, the above system is being described as a perturbed one obtained at uniform samples so that its fundamental matrix is time-varying (while constant in- between any two consecutive samples) so that the above results could be applied for the state and output trajectory solutions. A positive n- real matrix $M$ (respectively, n-real vector v) is denoted as $M \in \mathbf{R}_+^{n \times n}$, or simply as $M > 0$, (respectively as $v \in \mathbf{R}_+^n$, or simply as $v > 0$) what implies that all their entries are nonnegative and at least one is positive. The strongest notations $M \gg 0$, $v \gg 0$ stand when all the corresponding entries are positive. The notations $M \geq 0$, $v \geq 0$ stand for the matrix or vector being positive or null. This can be also denoted as $M \in \mathbf{R}_{0+}^{n \times n}$, $v \in \mathbf{R}_{0+}^{n \times n}$ (i.e. real matrix or vector with all its components being nonnegative). The subsequent related basic preliminary result follows by simple inspection of (10):

**Lemma 1**. Assume that $V(t) \in [0,1]$; $\forall t \in \mathbf{R}_{0+}$. Then for each fixed $t \in \mathbf{R}_{0+}$, one gets:

$A_1(x(t),t)$ and $A_4(x(t),t)$ are Metzler matrices if and only if $min(\omega, \beta, \sigma, \gamma(1-\rho)) \geq 0$

$A_2(x(t),t)$ and $A_3(x(t),t)$ are Metzler matrices if and only if $\beta = 0$ and $min(\omega, \sigma, \gamma(1-\rho)) \geq 0$

$A_1^0(x(t),t)$ and $A_4^0(x(t),t)$ are Metzler matrices if and only if $min(\nu, \omega+\nu, \beta, \sigma, \gamma(1-\rho)) \geq 0$

$A_2^0(x(t),t)$, $A_3^0(x(t),t)$ are Metzler matrices if and only if $min(\nu-\beta, \nu, \omega+\nu, \beta, \sigma, \gamma(1-\rho)) \geq 0$

$b_a \notin \mathbf{R}_+^4$; $c \in \mathbf{R}_+^4$ and c>>0; $\zeta > 0$ if $\nu > 0$ and $u(t) > 0$ if and $V(t) \in [0,1]$.

As a result, $0 \leq min(S(t), E(t)I(t), R(t)) \leq max(S(t), E(t)I(t), R(t)) \leq N(t)$; $\forall t \in \mathbf{R}_{0+}$. □

**Remark 1**. Note the following facts from Lemma 1 and the physical context of the SEIR – model:

a) All the real parameters which parameterize the SEIR-model (1)-(4) should be nonnegative. Therefore, the condition for $A_2(x(t),t)$ and $A_3(x(t),t)$ to be a Metzler matrix is not useful since the transmission constant β is required to be zero. If β = 0 then there is no coupling of the susceptible with the infected (see (1)-(2)) and the SEIR-model (1)-(4) is nonsense. Also, the condition $\gamma(1-\rho) \geq 0$ is useful only under $\gamma \geq 0$ and $\rho \in [0,1]$.

b) Note that $b_a \notin \mathbf{R}_+^4$, $\delta(t) := \zeta y(t) = \nu N(t) e_1 > 0$ if $\nu > 0$ and $V(t) \in [0,1]$, what implies that new defined control vector $u(t) = (1-V(t), 0, 0, V(t))^T$ is a nonnegative function. Then the positivity of the SEIR – model might be better characterized from (6) (instead of from (5) or(7)) and (8) in terms that the dynamic matrix is a Metzler time-invariant one and the control vector is positive as it is its coefficient by invoking the appropriate positivity conditions of Lemma 1. □

In the following, the matrices $A_4(x(t),t)$ and $A_4^0(x(t),t)$ are chosen to describe the above system. The subsequent study could be performed equivalently by using instead $A_1(x(t),t)$ and $A_1^0(x(t),t)$ Their associate subscript "4" is omitted in the notation as it is omitted the explicit dependence on x(t) of the matrix of dynamics, as follows $A(t) \equiv A_4(x(t), t)$ and $A^0(t) \equiv A_4^0(x(t), t)$ so that the SEIR-model (1)-(4) becomes compacted as:

$$\dot{x}(t) = A(t)x(t) + \nu N(t)u(t) \tag{11}$$

$$y(t) = c^T x(t) \tag{12}$$

where $u(t) = (1 - V(t), 0, 0, V(t))^T$ is the vaccination control generated for each given vaccination function $V(t)$. Also, note that

$$\dot{y}(t) = \dot{N}(t) = c^T \dot{x}(t) = c^T(A(t)x(t) + \nu y(t)u(t)) = (\nu - \mu)N(t) - \rho \gamma I(t) = c_1^T x(t) \tag{13}$$

where $c_1 = (\nu - \mu, \nu - \mu, \nu - \mu - \rho\gamma, \nu - \mu)^T$.

**Theorem 1**. Assume that $min(\beta, \sigma, \omega, \gamma) \geq 0$ and $\rho \in [0,1]$. Then, there exists a non-unique vaccination control $V \in PC^{(0)}(R_{0+}; [0,1])$ such that the state and output solution trajectories of the SEIR –model (1)-(4), equivalently (11)-(12), are nonnegative for all time for each given set of nonnegative initial conditions and any piecewise continuous vaccination function taking values in $[0,1]$ everywhere in its definition domain $R_{0+}$.

**Proof**: First, note that the second and third components of the state vector are the infected and infectious, respectively. It is clear that they are not directly affected by the vaccination function but only through the coupling with the susceptible related to the value of the disease transmission constant $\beta$. Note that if for both $i = 2,3$, $x_i(t) = 0 \Rightarrow \dot{x}_i(t) = 0 \Rightarrow x_i(t') = 0; \forall t' \geq t$ by simple inspection of (2)-(3). Note also that the state vector is continuous through time so that if only one of the populations of infected or infectious reaches zero at a time t by the first time then it follows also from (2)-(3) that its time - derivative is non less than zero at that time so that it cannot reach any negative value at $t^+$ (denoting the time argument for right limit of all partial population values at time t) under non-negative initial conditions. As a result, the infected and infectious cannot be negative at any time irrespective of the vaccination function if the susceptible population is non- negative at that time. Thus, if $S(t) \geq 0$; $\forall t \in R_{0+}$ for some vaccination function then $E(t) \geq 0$; $I(t) \geq 0$ for all time. Now, define the following matrices for any time $t \in [t_k, t_{k+1})$ for some ( in general time-varying) sampling period $T_k := t_{k+1} - t_k > 0$, where $SI := \{t_k\}_{k \in SZ_+}$ and $SP := \{T_k\}_{k \in SZ_+}$ are respectively the sequences of sampling instants and time periods where $SZ_+$ is some finite (and eventually empty if there is no sampling process) or infinite subset of $Z_+$ depending on the number of samples being finite or infinity. Define also the following real sequences at sampling instants $\{S_k\}_{k \in SZ_+}$, $\{E_k\}_{k \in SZ_+}$, $\{I_k\}_{k \in SZ_+}$, $\{R_k\}_{k \in SZ_+}$, $\{N_k\}_{k \in SZ_+}$:

$$A_k = A(t_k) := \begin{bmatrix} -\left(\mu + \frac{\beta I_k}{N_k}\right) & 0 & 0 & \omega \\ 0 & -(\mu+\sigma) & \frac{\beta S_k}{N_k} & 0 \\ 0 & \sigma & -(\mu+\gamma) & 0 \\ 0 & 0 & \gamma(1-\rho) & -(\mu+\omega) \end{bmatrix}$$

$$\tilde{A}_k(t) = A(t) - A_k := \begin{bmatrix} \beta\left(\frac{I_k}{N_k} - \frac{I(t)}{N(t)}\right) & 0 & 0 & 0 \\ 0 & 0 & -\beta\left(\frac{S_k}{N_k} - \frac{S(t)}{N(t)}\right) & 0 \\ 0 & 0 & 0 & 0 \\ 0 & 0 & 0 & 0 \end{bmatrix}$$

so that $A_k$ is a Metzler matrix ; $\forall k \in SZ_+$ so that $\tilde{A}_k(t) \in \mathbf{R}_{0+}^{n \times n}$ provided that $S_k \leq \frac{N_k S(t)}{N(t)}$ and $I_k > \frac{N_k I(t)}{N(t)}$ ; $\forall k \in SZ_+$ and $SI$ (and then $SZ_+$) is nonempty.

**Remark 2**. Note that in common cases $\frac{S(0)}{N(0)}$ and $\frac{S(\infty)}{N(\infty)}$ are unit or close to unit in the sense that all or almost all the population is susceptible when the disease propagation begins and this also happens at the end of the disease period when a new propagation cycle starts. There is also a finite time interval where the susceptible are decreasing towards numbers close to zero. After a minimum quotient value is reached there are time intervals where this quotient increases so that the above $S_k \leq \frac{N_k S(t)}{N(t)}$ assumption is feasible by choosing the elements of SI within such a time intervals.     □

Eq. (11) is equivalent to

$$\dot{x}_k(t) = A_k x_k + \left[v N_k(t) u_k(t) + \tilde{A}_k(t) x_k(t) + A_k(x_k(t) - x_k)\right] \quad (14)$$
$$= A_k(t) x_k(t) + v N_k(t) u_k(t)$$

with $x_k \equiv x(kT)$, $x_k(t) \equiv x(t)$ and $u_k(t) \equiv u(t)$; $\forall t \in [t_k, t_{k+1})$ where $SI \ni t_{k+1} = t_k + T_k$ is the next sampling instant to $t_k \in SI$ and $T_k \in SP$. Such an abuse of notation allows to refer signals defined for $t \in [t_k, t_{k+1})$ to the sampling instant $t_k$. Consider the positive real sequence of sampling periods $SP = \{T_k\}_{k \in SZ+}$ such that $T_k \in (0, \varepsilon_T)$; $\forall k \in SZ_+$. It follows from (14) that for any piecewise continuous vaccination function $V: \mathbf{R}_{0+} \to [0,1]$

$$\|x_k(t) - x(kT)\| \leq \varepsilon_T M_k\left(\|x_k(t)\|, \|u_k(t)\|, \varepsilon_{kT}\right) \quad (15)$$

where $M_k\left(\|x_k(t)\|, \|u_k(t)\|, \varepsilon_{kT}\right)$ is a norm-dependent monotone (in general non-strictly) continuous increasing real function of $\varepsilon_{kT}$ fulfilling $M_k\left(\|x_k(t)\|, \|u_k(t)\|, 0\right) = 0$ defined by

$$M_k(\|x_k(t)\|,\|u_k(t)\|,\varepsilon_{kT})=\max\|A_k x_k\|+2\nu \max_{\zeta\in[t_k,t_{k+1}]}\|x_k(\zeta)\|$$

$$+\max_{\zeta\in[t_k,t_{k+1}]}\|\widetilde{A}_k(\zeta)x(\zeta)\|+\max_{\zeta\in[t_k,t_{k+1}]}\|A_k(x_k(\zeta)-x_k)\|$$

$$\eta_k(t):=\left[\nu N_k(t)u_k(t)+\widetilde{A}_k(t)x_k-\left(\|A_k\|M_k(\|x_k(t)\|,\|u_k(t)\|,\varepsilon_{kT})+\|\widetilde{A}_k(t)\vartheta_{k1}(t)\|\right)e+\vartheta_k(t)\right]$$

for $\vartheta_{k1}(t):=x_k(t)-x_k$; $\forall k\in SZ_+$ and some $\vartheta_k(t)\in R_{0+}^n$. Since $\widetilde{A}_k(t)\in R_{0+}^{n\times n}$; $\forall k\in SZ_+$ and of norm being a continuous function of $\varepsilon_{kT}$ of zero value at $\varepsilon_{kT}=0$, i.e. for $t=kT$ and $M_k$ is a continuous monotone increasing function with $\varepsilon_{kT}$ being zero for $\varepsilon_T=0$, it follows that $\eta_k(t)\in R_{0+}^4$; $\forall t\in[t_k,t_{k+1}]$ for $\varepsilon_{kT}\in[0,\bar{\varepsilon}_k)$; $\forall k\in SZ_+$ with a sufficiently small $\bar{\varepsilon}_k\in R_+$ so that from (14)-(15):

$$\dot{x}_k(t)=A_k x_k+\eta_k(t)$$

$$\Rightarrow \left(x_k>0\Rightarrow x_k(t)=e^{A_k(t-kT)}x_k+\int_0^{t-kT}e^{A_k(t-kT-\tau)}\eta_k(kT+\tau)d\tau\right); \forall t\in[t_k,t_{k+1})$$

; $\forall t_k\in SI$ for some vaccination function taking values in $[0,1]$ which keeps nonnegative also the susceptible and immune for all time apart from the previously proven nonnegative of infected and infectious. By using complete induction for the next interval provided that the property holds for all the previous ones, the positivity is proven for all time. Since $y(t)=N(t)=c^T x(t)$ then $x(t)\geq 0\Rightarrow y(t)\geq 0$ for any $t\in R_{0+}$ since $c>>0$. □

A simpler alternative proof of Theorem 1 is now provided by decomposing $A(t)\equiv A^*+B(t)$ in a non-unique way, where

$$A^*:=\begin{bmatrix}-\left(\mu+\dfrac{\beta I^0}{N^0}\right) & 0 & 0 & \omega \\ 0 & -(\mu+\sigma) & 0 & 0 \\ 0 & \sigma & -(\mu+\gamma) & 0 \\ 0 & 0 & \gamma(1-\rho) & -(\mu+\omega)\end{bmatrix}; B(t):=\begin{bmatrix}\dfrac{\beta I^0}{N^0}-\dfrac{\beta I(t)}{N(t)} & 0 & 0 & 0 \\ 0 & 0 & \dfrac{\beta S(t)}{N(t)} & 0 \\ 0 & 0 & 0 & 0 \\ 0 & 0 & 0 & 0\end{bmatrix}$$

(16)

so that the system (11)-(12) can be rewritten equivalently as:

$$\dot{x}(t)=A^* x(t)+B(t)x(t)+\nu y(t)u(t) \quad (17)$$

$$y(t)=c^T x(t) \quad (18)$$

where $u(t)=(1-V(t),0,0,V(t))^T$ and $y(t)=N(t)$ with $I^0$ and $N^0$ being positive constants being arbitrary except for the constraint $N^0\geq I^0\geq \dfrac{N^0}{N(t)}I(t)$ so that $B(t)\geq 0$. For instance, $N^0$ may be chosen to be a known upper-bound of $\max_{t\in R_{0+}}N(t)$, for instance, $N_0=N(0)$ in the common case that $\dot{N}(t)\leq 0$; $\forall t\in R_{0+}$.

**Alternative Proof of Theorem 1**: The unique state/output trajectory solutions of the dynamic system (17)-(18), equivalent to (11)-(12), for each set of initial conditions and each vaccination function are given by:

$$x(t) = e^{A^* t}\left(x_0 + \int_0^t e^{-A^* \tau}(B(\tau)x(\tau) + \nu N(\tau)u(\tau))d\tau\right) ; \forall t \in \mathbf{R}_{0+} \quad (19)$$

$$y(t) = c^T e^{A^* t}\left(x_0 + \int_0^t e^{-A^* \tau}(B(\tau)x(\tau) + \nu N(\tau)u(\tau))d\tau\right) ; \forall t \in \mathbf{R}_{0+} \quad (20)$$

so that

$$[x_0 \geq 0 \wedge (V(t) \in [0,1] ; \forall t \in \mathbf{R}_{0+})] \Rightarrow [u(t) \geq 0 \wedge (y(t) \equiv N(t) \geq 0) \wedge (x(t) \geq 0; \forall t \in \mathbf{R}_{0+})] \quad (21)$$

since $A^*$ is a constant Metzler matrix by construction so that the fundamental matrix $\Psi(t) := e^{A^* t}$ is positive, see (16), i.e. $\Psi : \mathbf{R}_{0+} \to \mathbf{R}_+^{4 \times 4}$ (or $\Psi(t) > 0$; $\forall t \in \mathbf{R}_{0+}$), $\nu \geq 0$, $\beta \geq 0$, $c \in \mathbf{R}_+^4$ (with $c \gg 0$), $B : \mathbf{R}_{0+} \to \mathbf{R}_{0+}^4$ (i.e. $B(t) \geq 0$, $\forall t \in \mathbf{R}_{0+}$ see (16)). Since (19) and (20) are unique for each set of nonnegative initial conditions and vaccination function in $[0,1]$, the positivity properties are independent of the decomposition (16) and then valid for (11)-(12) as a result. □

## 4. Stability

Two simple stability results follow which are based on the conditions for joint achievement of the positivity property and the uniform boundedness for all time of the total population N(t).

**Theorem 2**. Assume that $0 \leq \nu \leq \mu$, $\min(\beta, \sigma, \omega, \gamma) \geq 0$, $\rho \in [0,1]$, $V \in PC^{(0)}(\mathbf{R}_{0+}; [0,1])$ and $\frac{I^0}{N^0} \geq \frac{I(t)}{N(t)}$; $\forall t \in \mathbf{R}_{0+}$. Then, the SEIR- model (1)-(4) is positive, stable and $N(t)$ converges asymptotically to a finite limit $N(\infty) \leq N(0)$ as $t \to \infty$. If $\mu = \nu$, and also $\rho = 0$ or $\gamma = 0$ then $N(t) = N(0); \forall t \in \mathbf{R}_{0+}$.

**Proof**: The above conditions guarantee that the model is positive (see Theorem 1 and Lemma 1) guaranteeing that $A^*$ is a Metzler matrix and $B(t) \geq 0$; $\forall t \in \mathbf{R}_{0+}$). Thus, $\gamma \geq \rho \gamma \geq 0$ and also no population (susceptible, infected, infectious and immune) may be negative at any time. As a result, all those populations are uniformly bounded for all time if $N(t)$ is uniformly bounded for all time. Thus, since $\dot{y}(t) = \dot{N}(t) = (\nu - \mu)N(t) - \rho \gamma I(t)$ and since the system is positive from Theorem 1, one gets:

$$[(0 \leq \nu \leq \mu) \wedge (\rho \gamma \geq 0)]$$
$$\Rightarrow \dot{N}(t) = (\nu - \mu)N(t) - \rho \gamma I(t) \leq 0 \Rightarrow 0 \leq N(t) \leq N(0) < \infty ; \forall t \in \mathbf{R}_{0+}$$
$$\Rightarrow 0 \leq \max(S(t), E(t), I(t), R(t)) \leq N(0) < \infty$$

Since N(t) is monotone decreasing on its definition domain $\mathbf{R}_{0+}$ then it converges to a limit not higher than $N(0)$. If $\mu - \nu = \rho \gamma = 0$ then $\dot{N}(t) \equiv 0 \Rightarrow N(t) = N(0); \forall t \in \mathbf{R}_{0+}$. □

**Theorem 3**. Assume that $\nu > \mu$, $min(\beta, \sigma, \omega, \gamma) \geq 0$, $\rho \in [0,1]$ and $V \in PC^{(0)}(\mathbf{R}_{0+}; [0,1])$. Then, the following properties hold:

(i) The SEIR- model (1)-(4) is positive.

(ii) A necessary condition for the stability of the SEIR model (1)-(4) is $\gamma \geq \dfrac{\nu - \mu}{\rho}$ subject to $\rho > 0$ including the case that the total population N(t) extinguishes in finite time. As a result, the set $\Gamma := \left\{ t \in \mathbf{R}_{0+} : I(t) < \dfrac{(\nu-\mu)N(t)}{\rho\gamma} \right\} \subseteq \mathbf{R}_{0+}$ is either empty or connected of finite Lebesgue measure or non-connected with each connected component being of finite measure.

(iii) The SEIR – model (1)-(4) is stable if and only if
$$N(0) = \rho\gamma \int_0^\infty e^{(\mu-\nu)\tau} I(\tau)d\tau = \rho\gamma \int_0^\infty e^{(\mu-\nu)\tau} e_3^T x(\tau)d\tau$$

**Proof**: Property (i) follows directly from Theorem 1. To prove (ii), note that the uniform boundedness of $N(t)$ for all time requires that $\dot{N}(t) = (\nu - \mu)N(t) - \rho\gamma I(t) \leq 0$; $\forall t \in T_s \subseteq \mathbf{R}_{0+}$ with $T_s \equiv \bigcup_{i \in N_s} T_{is}$ being non-empty and subject to $\mathsf{M}(T_s) = \infty$, i.e. it is of infinite Lebesgue measure where the disjoint connected components $T_{is}$ of $T_s$ can be of finite or infinite measure but if its number $Z_s \subseteq \mathbf{Z}_+$ is finite, i.e. $1 \leq c := card\, Z_s < \infty$ then $\mathsf{M}(T_{cs}) = \infty$. This follows by contradiction. Assume that $c = 0$ then $\dot{N}(t) = (\nu - \mu)N(t) - \rho\gamma I(t) > 0$; $\forall t \in \mathbf{R}_{0+}$. Then, $N(t) \to \infty$ as $t \to \infty$ and the SEIR- model is unstable. Now, assume that $T_{cs} = [t_{1cs}, t_{2cs}]$ with $1 < c := card\, Z_s < \infty$ and $\mathsf{M}(T_{cs}) = |t_{2cs} - t_{1cs}| < \infty$. Then, $\dot{N}(t) = (\nu - \mu)N(t) - \rho\gamma I(t) > 0$; $\forall t \geq t_{c2}$ so that again $N(t) \to \infty$ as $t \to \infty$. As a result, the condition below is required for stability

$\dot{N}(t) = (\nu - \mu)N(t) - \rho\gamma I(t) \leq 0$; $\forall t \in T_s$ with $\mathsf{M}(T_s) = \infty$ what implies

$I(t) \geq \dfrac{(\nu-\mu)N(t)}{\rho\gamma} = \dfrac{(\nu-\mu)(S(t)+E(t)+I(t)+R(t))}{\rho\gamma}$; $\forall t \in T_s \Leftrightarrow I(t) \geq \dfrac{(\nu-\mu)(S(t)+E(t)+R(t))}{\rho\gamma+\mu-\nu}$; $\forall t \in T_s$

so that, since the model is positive, either $\gamma \geq \dfrac{\nu - \mu}{\rho}$ or $S(t) + E(t) + R(t) = 0$ so that $I(t) = N(t)$; $\forall t \in T_s$. In the first case, Property (ii) is proven. The second one is impossible provided that Property (iii), what is being proven next, is true except if $I(t) = N(t) \equiv 0$ on a set of infinite measure. To prove Property (iii), calculate the time-integral of $\dot{N}(t) = (\nu - \mu)N(t) - \rho\gamma I(t)$ to yield:

$$e^{-(\nu-\mu)t} N(t) = \left( N(0) - \rho\gamma \int_0^t e^{-(\nu-\mu)\tau} I(\tau)d\tau \right); \forall t \in \mathbf{R}_{0+} \qquad (22)$$

Assume that the SEIR- model (1)-(4) is stable. Then, $N(t)$ is uniformly bounded for all time and all the partial populations of susceptible, infected, infectious and immune are nonnegative and uniformly bounded for all time since the SEIR – model is positive. Taking limits in both sides of (22) as $t \to \infty$ yields $N(0) = \rho\gamma \int_0^\infty e^{-(v-\mu)\tau} I(\tau) d\tau$ since $v > \mu$. The sufficiency part of the property has been proven. The necessity part is proven by contradiction. Assume that $N(0) < \rho\gamma \int_0^\infty e^{-(v-\mu)\tau} I(\tau) d\tau$. Then, $N(t) < 0$ for sufficiently large finite time what contradicts the positivity of the SEIR -model. Now, assume that $\infty > N(0) > \rho\gamma \int_0^\infty e^{-(v-\mu)\tau} I(\tau) d\tau$ then $N(t) \to \infty$ as $t \to \infty$ so that the model is unstable. As a result, it is necessary for stability that $N(0) = \rho\gamma \int_0^\infty e^{-(v-\mu)\tau} I(\tau) d\tau$. Property (iii) has been proven and then the remaining reasoning to complete the proof of Property (ii). □

The subsequent result is related to stability irrespective of the positivity of the SEIR – model. Therefore, it is not directly applicable to real problems.

**Theorem 4**. Assume that $min(\beta,\sigma,\omega,\gamma) \geq 0$, $\rho \in [0,1]$ and $V \in PC^{(0)}(\boldsymbol{R}_{0+}; [0,1])$. Assume also that $\mu > 0$, $max(\sigma,\gamma) < -\mu$ and, furthermore, that either (1) $0 \leq v \leq \mu$, or (2) $\mu > 4\beta + v$, $v > \beta \geq 0$. Assume also that the vaccination function is piecewise continuous taking values in $[0,1]$ everywhere in $\boldsymbol{R}_{0+}$. Then, the SEIR – model (1)-(4) is stable.

**Proof**: The SEIR model is positive according to Theorem 1. Since $\dot{y}(t) = \dot{N}(t) = (v-\mu)N(t) - \rho\gamma I(t)$ then since the system is positive from Theorem 1, on gets:

$[(0 \leq v \leq \mu) \wedge (\rho\gamma \geq 0)]$
$\Rightarrow \dot{N}(t) = (v-\mu)N(t) - \rho\gamma I(t) \leq 0 \Rightarrow 0 \leq N(t) \leq N(0) < \infty \Rightarrow 0 \leq max(S(t),E(t),I(t),R(t)) \leq N(0) < \infty$

$\forall t \in \boldsymbol{R}_{0+}$ so that the SEIR model is stable. Now, assume that $v > \mu > v - \beta I^0/N^0 > 0$, the lower-bound being always feasible by appropriately fixing $N^0 \geq I^0 \geq N^0 \left( \underset{t \in \boldsymbol{R}_{0+}}{max} \frac{I(t)}{N(t)} \right)$. In this case, $B(t) \geq 0$ requires $\beta I(t)/N(t) < v$ which is guaranteed if $\beta < v$. Note that Eq. (17) may be rewritten equivalently as follows:

$$\dot{x}(t) = \left( A^{0*} + B(t) \right) x(t) + b_a y(t) V(t) \tag{23}$$

where $A^{0*} = A^* + \zeta c^T$ is a constant real stability matrix of eigenvalues $-\left( \mu + \frac{\beta I^0}{N^0} - v \right)$, $-(\mu+\sigma)$, $-(\mu+\gamma)$, $-(\mu+\omega)$. Since the system is positive, all the populations in (1)-(4) are nonnegative and less than N(t) for all time. Then, $0 \leq I(t)/N(t) \leq 1$ for all time and B(t) in (16) is uniformly bounded on $R_{0+}$ as a result. The unforced system of (23) is $\dot{z}(t) = (A^{0*} + B(t))z(t)$. The norm of the evolution operator $\Phi(t,0)$ defining the state –trajectory solution of the unforced system is uniformly bounded for all time by $K_\Phi e^{-(\rho_0 - b)t}$ where $K_\Phi \geq 1$ is norm-dependent, $b > 0$ is a real constant related to some maximum norm $\max_{t \in R_{0+}} \|B(t)\|$, and $\rho_0 > 0$ is the minus stability abscissa of the stability matrix $A^{0*}$; i.e. any constant less than its minus minimum eigenvalue (or equal to it if its multiplicity is one). For instance, by using the $\ell_2$- (spectral) matrix norm and the property $b := \max_{t \in R_{0+}} \|B(t)\|_2 \leq 2\|B(t)\|_\infty$ (since the order of B(t) is n=2), one gets $b \leq 4\beta$ and $\rho_0 > b$ if $\mu > 4\beta + \nu$ what implies global asymptotic Lyapunov´s stability of the unforced system of (23). The forced system is then globally Lyapunov´s stable since the forcing function is uniformly bounded for all time. □

## 5. Control via vaccination

The vaccination problem is now focused on. *The basic objective is the achievement either asymptotically or in finite time that the whole population be immune* irrespective of the initial conditions. It is assumed that the whole population N (t) and the infectious one I (t) are known through time what is a reasonable proposition. The proposed vaccination law has the following structure:

$$V(t) = \begin{cases} V_a(t) & \text{if } V_a(t) \in [0,1] \\ 1 & \text{if } V_a(t) > 1 \\ 0 & \text{if } V_a(t) < 0 \end{cases} \quad (24)$$

$$V_a(t) = \frac{1}{\nu N(t)} \left( K_N(t)N(t) + K_I(t)I(t) + K_R(t)R^*(t) + K_{Rd}(t)\dot{R}^*(t) \right) \quad (25)$$

$$R^*(t) = h(t)N(t) \; ; \; R^*(0) = R(0) \quad (26)$$

where $R^*(t)$ is the desired reference immune population and $h \in PC^{(1)}(R_{0+}; R_{0+})$, namely, it is time- differentiable with piecewise continuous time derivative everywhere in its definition domain, which is monotone increasing (including a potential constant choice) and converging to one in order to achieve $R^*(t) \to N(t)$ as $t \to \infty$. The controller is defined by time-varying controller real gains $K_N(t)$, $K_I(t), K_R(t)$ and $K_{Rd}(t)$, which satisfy the following design constraints:

$$K_N(t) = -\left(K_R(t) + (\nu - \mu)K_{Rd}(t)\right)h(t) - K_{Rd}(t)\dot{h}(t) + \varepsilon_0(1 - \varepsilon g(t)) \quad (27.a)$$

$$K_I(t) = \gamma \rho K_{Rd}(t)h(t) \quad (27.b)$$

$$\dot{h}(t) = \frac{1}{K_{Rd}(t)} \left[ K_R(t)(R^*(t) - h(t)N(t)) + K_{Rd}(t)(\dot{R}^*(t) + (\mu - \nu)N(t)h(t)) \right]$$

$$= \frac{1}{K_{Rd}(t)} \left[ K_R(t)R^*(t) + K_{Rd}(t)\dot{R}^*(t) + ((\mu - \nu)K_{Rd}(t) - K_R(t))N(t)h(t) \right] \quad (28)$$

where (13) has been used and $\varepsilon \in R_{0+}$, $\varepsilon_0 \in R_+$ and $g : R_{0+} \to R$ if $\varepsilon > 0$ and $g \equiv 0$ if $\varepsilon = 0$ is a design function which modulates the desired rate to reach the objective of asymptotically tracking the whole population by the immune one. The substitution of (27)-(28) into (24)-(26) yields:

$$v N(t)V_a(t) = \varepsilon_0 (1 - \varepsilon g(t))N(t) \qquad (29)$$

$$v N(t)V(t) = \varepsilon_0 (1 - \varepsilon g(t))N(t) + v N(t)[\vartheta_1(t) - (\vartheta_1(t) + \vartheta_0(t))V_a(t)] \qquad (30)$$

where the indicator functions are defined as follows for each $t \in R_{0+}$:

$\vartheta_1(t) = 1$ if $V_a(t) > 1$ and $\vartheta_1(t) = 0$ if $V_a(t) \leq 1$

$\vartheta_0(t) = 1$ if $V_a(t) < 0$ and $\vartheta_0(t) = 0$ if $V_a(t) \geq 0$

The substitution of (29)-(30) in (4) yields the following controlled evolution of the immune population:

$$\begin{aligned}\dot{R}(t) &= -(\mu + \omega)R(t) + \varepsilon_0 (1 - \varepsilon g(t))N(t) + v N(t)[\vartheta_1(t) - (\vartheta_1(t) + \vartheta_0(t))V_a(t)] + \gamma (1 - \rho)I(t) \\ &= -(\mu + \omega)R(t) + \varepsilon_0 (1 - \varepsilon g(t))N(t) + v N(t)\vartheta_1(t) \\ &\quad + \gamma (1 - \rho)I(t) - \varepsilon_0 (\vartheta_1(t) + \vartheta_0(t))(1 - \varepsilon g(t))N(t) \\ &= -(\mu + \omega)R(t) + \varepsilon_0 (1 - \varepsilon g(t))N(t) + [(v - \varepsilon_0) + \varepsilon_0 \varepsilon g(t)]N(t)\vartheta_1(t) \\ &\quad - \varepsilon_0 (1 - \varepsilon g(t))N(t)\vartheta_0(t) + \gamma (1 - \rho)I(t) \end{aligned} \qquad (31)$$

Eq. 31 has the following particular cases which will be shown to be of design interest when performing particular vaccination strategies:

**1)** $\dot{R}(t) = -(\mu + \omega)R(t) + \varepsilon_0 (1 - \varepsilon g(t))N(t)$ \qquad (32)

which holds by comparison with (31) if the following identity holds:

$$(v - \varepsilon_0)N(t)\vartheta_1(t) + \varepsilon_0 \varepsilon g(t)N(t)\vartheta_1(t) - \varepsilon_0 N(t)(1 - \varepsilon g(t))\vartheta_0(t) + \gamma (1 - \rho)I(t) = 0 ; \forall t \in R_{0+}$$

which is achieved through the subsequent choice of the control function $g : R_{0+} \to [0, \varepsilon^{-1}]$:

$$g(t) = \frac{(\varepsilon_0 \vartheta_0(t) + (\varepsilon_0 - v)\vartheta_1(t))N(t) - \gamma (1 - \rho)I(t)}{\varepsilon_0 \varepsilon N(t)(\vartheta_0(t) + \vartheta_1(t))} \text{ if } \vartheta_0(t) + \vartheta_1(t) = 1 ; \text{ i.e. if } V(t) \neq V_a(t)$$

; $\forall t \in R_{0+}$ \qquad (33.a)

$$g(t) = \frac{1}{\varepsilon}\left(1 - \frac{\gamma (1 - \rho)I(t)}{\varepsilon_0 N(t)}\right) \text{ if } \vartheta_0(t) = \vartheta_1(t) = 0 ; \text{ i.e. if } V(t) = V_a(t) ; \forall t \in R_{0+} \qquad (33.b)$$

provided that $\varepsilon_0 > max(v, \gamma (1 - \rho))$. Note that the above constraint replaced in the expression of the controller gain $K_N(t)$ Eq. (27) has the suitable advantage that the auxiliary control $V_a$ equalizes the true one taking values in $[0,1]$ so that both stability and positivity of the SEIR- model are guaranteed and the controller is potentially useful to asymptotically achieve total immunity of the whole population.

2) $\dot{R}(t) = -(\mu+\omega)R(t) + \nu N(t)\vartheta_1(t) + \gamma(1-\rho)I(t)$ (34)

$\leq -(\mu+\omega)R(t) + (\nu + \gamma(1-\rho))N(t)$

if $g(t) = g \equiv 1/\varepsilon$ ; $\forall t \in R_{0+}$ ( i.e. it is is constant). If, in addition, $\vartheta_1 \equiv 0$, i.e. $V_a(t) \in (-\infty, 1]$ ; $\forall t \in R_{0+}$ then

$\dot{R}(t) = -(\mu+\omega)R(t) + \gamma(1-\rho)I(t) \leq -(\mu+\omega)R(t) + \gamma(1-\rho)N(t)$ (35)

3) $\dot{R}(t) = -(\mu+\omega)R(t) + \nu N(t) + \nu N(t)(\varepsilon g(t)(\vartheta_1(t)-1) - (1-\varepsilon g(t))\vartheta_0(t)) + \gamma(1-\rho)I(t)$

(36)

$\leq -(\mu+\omega)R(t) + \nu N(t) + \nu(\varepsilon g(t) + \gamma(1-\rho))N(t)$

if $\nu = \varepsilon_0$. New particular cases are obtained from (34) as follows:

a) If, in addition, $\vartheta_0 \equiv 0, \vartheta_1 \equiv 0$; i.e. $V \equiv V_a \in [0,1]$ then the choice $g(t) = \dfrac{\gamma(1-\rho)I(t)}{\varepsilon \nu N(t)}$ yields

$\dot{R}(t) = -(\mu+\omega)R(t) + \nu N(t)$ (37)

b) If $\vartheta_0 \equiv 1, \vartheta_1 \equiv 0$; i.e. $V_a(t) < 0; \forall t \in R_{0+}$ then the choice $g(t) = \dfrac{1}{\varepsilon}\left(\dfrac{\gamma(1-\rho)I(t)}{\nu N(t)} - 1\right)$ yields

$\dot{R}(t) = -(\mu+\omega)R(t) + \gamma(1-\rho)I(t)$ (38)

c) If $\vartheta_0 \equiv 0, \vartheta_1 \equiv 1$; i.e. $V_a(t) > 1 \, \forall t \in R_{0+}$ then (36) is obtained irrespective of the function g.

d) If $g \equiv 0$ then

$\dot{R}(t) = -(\mu+\omega)R(t) + \nu N(t)(1-\vartheta_0(t)) + \gamma(1-\rho)I(t) \leq -(\mu+\omega)R(t) + (\nu + \gamma(1-\rho))N(t)$ (39)

e) $\dot{R}(t) = -(\mu+\omega)R(t) + \gamma(1-\rho)I(t) \leq -(\mu+\omega)R(t) + \gamma(1-\rho)N(t)$ (40)

provided that $g \equiv 0$ , $\vartheta_0 \equiv 1$ ; i. e. $V_a(t) < 0$; $\forall t \in R_{0+}$.

The model could be useful for asymptotic immunity tracking of the whole population but its positivity for all time is not " a priori" guaranteed. The solution of the controlled immune population for any set of given initial conditions is calculated from (31) as follows:

$R(t) = e^{-(\mu+\omega)t}\bigg( R(0)$

$+ \int_0^t e^{(\mu+\omega)\tau}\Big\{ \varepsilon_0(1-\varepsilon g(\tau)) + [(\nu-\varepsilon_0) + \varepsilon_0 \varepsilon g(\tau)]\vartheta_1(\tau) - \varepsilon_0(1-\varepsilon g(\tau))\vartheta_0(\tau)$

$$+ \gamma(1-\rho)I(\tau)\} N(\tau)d\tau \Bigg) \tag{41}$$

The following assumption is coherent with real situations and useful for further study. It is concerned with the uniform boundedness for all time of the whole population which is also nonnegative. In other words, the disease can make the population to decrease or even to extinguish in the worst case but never to diverge.

**Assumption 1**. There exist real constants $N_i = N_i(N(0))$; i=1,2 subject to $0 \leq N_1 \leq N_2 < \infty$ such that $N(t) \in [N_1, N_2]$; $\forall t \in \mathbf{R}_{0+}$. □

It is direct the fact that Assumption 1 holds if Theorem 2 holds or if Theorem 3 (iii) holds. Theorem 3(ii) implies the following result:

**Assertion 1**. If Assumption 1 holds then either $0 \leq \nu \leq \mu$ or if $\nu > \mu \geq 0$ then $\rho > 0$ and, from Theorem 3, either $\frac{I(t)}{N(t)} \geq \frac{\nu - \mu}{\rho \gamma}$; $\forall t \geq t^*$ (some finite $t^* \in \mathbf{R}_{0+}$) or $\frac{I(t)}{N(t)} = \frac{\nu - \mu}{\rho \gamma}$; $\forall t \in S := \{t_i\}_{t_i \in \mathbf{R}_{0+}}$ with $|t_{i+1} - t_i| < \infty$, and

$$\sum_{\ell=0}^{i} \prod_{j=\ell+1}^{i} e^{(\nu-\mu)(t_{j+1}-t_j)} \left[ \int_0^{T_\ell} e^{(\nu-\mu)(T_\ell - \tau)} I(t_\ell + \tau)d\tau \right] \in \left[ \frac{e^{(\nu-\mu)t_{i+1}} N(0) - N_2}{\rho \gamma}, \frac{e^{(\nu-\mu)t_{i+1}} N(0) - N_1}{\rho \gamma} \right]$$
(42)

**Proof**: It follows from Theorem 3 (ii) since that property holds if $0 \leq \nu \leq \mu$; $\forall t \in \mathbf{R}_{0+}$. If $\nu > \mu \geq 0$ and $\frac{I(t)}{N(t)} \geq \frac{\nu - \mu}{\rho \gamma} \Leftrightarrow \dot{N}(t) \leq 0$; $\forall t \geq t^*$, then the uniform boundedness of N(t) is guaranteed for all time. Such a boundedness is also achievable if $\dot{N}(t)$ alternates its sign over consecutive time intervals of appropriate minimum size as follows. In that case, $I(t)$ has to take sufficiently large values at certain time intervals to generate a negative time-derivative to compensate former positive values. Then, it exists an infinite sequence of time instants $t_i$, $\frac{I(t_i)}{N(t_i)} = \frac{\nu - \mu}{\rho \gamma}$ implying that $\dot{N}(t_i) = 0$ while $\dot{N}(t)$ it has opposite signs in the former and next time interval $[t_{i-1}, t_i)$ and $[t_i, t_{i+1})$. By defining the time interval sequence $\{T_i := t_{i+1} - t_i\}_{t_i \in S}$, the constraint (42) follows from integrating (13) for all time to yield:

$$N(t_{i+1}) = e^{(\nu-\mu)t_{i+1}} N(0) - \rho \gamma \sum_{\ell=0}^{i} \prod_{j=\ell+1}^{i} e^{(\nu-\mu)(t_{j+1}-t_{j+1})} \left[ \int_0^{T_\ell} e^{(\nu-\mu)(T_\ell - \tau)} I(t_\ell + \tau)d\tau \right]; \forall t_i \in S$$

□

The asymptotic tracking property of the suitable reference of the immune population for the given vaccination law is formalized in the subsequent result:

**Theorem 5**. Under Assumption 1, the following particular cases of the vaccination control law (24) – (28):

(i) $\bar{R} := \dfrac{\varepsilon_0}{\mu+\omega}\left(1-\varepsilon \min_{t\in R_{0+}} g(t)\right) N_2$, subject to $\dfrac{\varepsilon_0}{\mu+\omega}\left(1-\varepsilon \min_{t\in R_{0+}} g(t)\right) \leq 1$, if g(t) is defined by

(33.a)

(ii) $\bar{R} := \dfrac{v+\gamma(1-\rho)}{\mu+\omega} N_2$, subject to $\dfrac{v+\gamma(1-\rho)}{\mu+\omega} \leq 1$, if $g(t) = \dfrac{1}{\varepsilon}$

(iii) $\bar{R} := \dfrac{\gamma(1-\rho)N_2}{\mu+\omega}$, subject to $\dfrac{\gamma(1-\rho)}{\mu+\omega} \leq 1$, if $g(t) = \dfrac{1}{\varepsilon}$ and $\vartheta_1 \equiv 0$

(iv) $\bar{R} := \dfrac{\gamma(1-\rho)N_2}{\mu+\omega}$, subject to $\dfrac{\gamma(1-\rho)}{\mu+\omega} \leq 1$, if $v = \varepsilon_0$, $\vartheta_0 = \vartheta_1 \equiv 0$ and $g(t) = 1/\varepsilon$

(v) $R(\infty) = 0$ if $v = \varepsilon_0$, $\gamma(1-\rho) = 0$ and $\vartheta_0 = 1; \vartheta_1 \equiv 0$. Furthermore, the immune population converges exponentially to zero irrespective of the initial conditions. If $\gamma(1-\rho) \neq 0$ then $\bar{R} := \dfrac{\gamma(1-\rho)N_2}{\mu+\omega}$.

(vi) $\bar{R} := \dfrac{v+\gamma(1-\rho)}{\mu+\omega} N_2$, subject to $\dfrac{v+\gamma(1-\rho)}{\mu+\omega} \leq 1$, if $v = \varepsilon_0$ and $g = \vartheta_0 \equiv 0$

(vii) $\bar{R} := \dfrac{\gamma(1-\rho)N_2}{\mu+\omega}$, subject to $\dfrac{\gamma(1-\rho)}{\mu+\omega} \leq 1$, if $v = \varepsilon_0$, $g \equiv 0$ and $\vartheta_0 \equiv 1$

(viii) $\bar{R} := \dfrac{(v+\gamma(1-\rho))+\varepsilon v \max_{t\in R_{0+}} g(t)}{\mu+\omega} N_2$, subject to $\dfrac{(v+\gamma(1-\rho))+\varepsilon v \max_{t\in R_{0+}} g(t)}{\mu+\omega} \leq 1$, if $v = \varepsilon_0$

yield the asymptotic upper-bounding tracking property of the immune $R(t) \to R(\infty) \leq \bar{R}$ as $t \to \infty$ at an exponential rate.

**Proof**: Property **(i)** follows from (32) –(33). Properties **(ii)-(viii)** follow, respectively from (34), (35), (37), (38), (39), (40) and (36). □

The following three results are related to the asymptotic convergence to zero of the immune population, its convergence to a finitely bounded limit with known upper and lower- bounds, and its convergence to the whole population provided to be constant, respectively.

**Theorem 6**. Assume that the vaccination control law (24)-(30) is generated with

$$g(\tau) = \dfrac{N(\tau) - e^{-\vartheta \tau}}{\varepsilon N(\tau)}; \quad \forall \tau \in [0, t); \quad \forall t \in R_{0+} \tag{43}$$

for some design real constant $\vartheta > \mu+\omega$ with $v \geq \varepsilon_0 \geq \dfrac{v N_1}{2N_1 - N_2} \geq v$ and $g(t) \leq \dfrac{\varepsilon_0 - v}{\varepsilon \varepsilon_0} \leq 0$;

$\forall t \in R_{0+}$. Then, the following two propositions hold:

**(i)** $R(t)$ converges exponentially to zero as $t \to \infty$ for any initial conditions of the SEIR – model (1)-(4).

**(ii)** Assume in addition that the total population $N(t)$ satisfies Assumption 1 for all time with $N_1 \leq N_2 < 2N_1$ and assume also that $h(t) = \dfrac{1}{N(t)}\left(e^{-(\mu+\omega)t}R(0) + \varepsilon_0 \dfrac{e^{-(\mu+\omega)t} - e^{-\vartheta t}}{\vartheta - \mu - \omega}\right)$

**Proof**: Note that $g(\tau) = \dfrac{N(\tau) - e^{-\vartheta\tau}}{\varepsilon N(\tau)} \leq \dfrac{N_2 - e^{-\vartheta\tau}}{\varepsilon N_1} \leq \dfrac{2\varepsilon_0 - \nu}{\varepsilon \varepsilon_0}$; $\forall \tau \in [0, t)$; $\forall t \in \mathbf{R}_{0+}$ since $\varepsilon_0 \geq \dfrac{\nu N_1}{2N_1 - N_2}$ implies that $V \equiv V_a \in [0, 1]$ so that the controlled SEIR-model via the given vaccination law is positive. One gets from the chosen vaccination control law (24)-(28) subject to (32) that:

$$R(t) = e^{-(\mu+\omega)t}\left(R(0) + \varepsilon_0 \int_0^t e^{(\mu+\omega)\tau}(1 - \varepsilon g(\tau))N(\tau)d\tau\right)$$

$$= e^{-(\mu+\omega)t}\left(R(0) + \varepsilon_0 \int_0^t e^{(\mu+\omega-\vartheta)\tau}d\tau\right) = e^{-(\mu+\omega)t}\left(R(0) + \varepsilon_0 \dfrac{1 - e^{-(\vartheta-\mu-\omega)t}}{\vartheta - \mu - \omega}\right)$$

Then $R(t) \to 0$ as $t \to \infty$ at exponential rate. □

**Corollary 1**. Assume that the immune population reference to be tracked by the vaccination law is

$$R^*(t) = h(t)N(t) = \varepsilon_0 \lim_{t \to \infty}\int_0^t e^{-(\mu+\omega)(t-\tau)}(1 - \varepsilon g(\tau))N(\tau)d\tau; \quad \forall t \in \mathbf{R}_{0+} \quad (44)$$

under the precise definition of the modulating function $h \in PC^{(1)}(\mathbf{R}_{0+}; \mathbf{R}_{0+})$. Then,

$R \equiv R^*$ implying also

$$R(t) \to R^*(t) = h(\infty)N(t) = \dfrac{\varepsilon_0 N(t)}{\vartheta - \mu - \omega} \in \left[\dfrac{\varepsilon_0 N_1}{\vartheta - \mu - \omega}, \dfrac{\varepsilon_0 N_2}{\vartheta - \mu - \omega}\right] \text{ as } t \to \infty$$

so that perfect tracking and asymptotic tracking of the immunity population are both achieved even if the total population N(t) does not have a limit. In the case that $\rho = \nu - \mu = 0$, $N(t) = N(0) = N = N_1 = N_2$; $\forall t \in \mathbf{R}_{0+}$, it follows that

$$R(t) = R^*(t) = h(t)N = R(\infty) = h(\infty)N(\infty) = \dfrac{\varepsilon_0 N}{\vartheta - \mu - \omega}; \quad \forall t \in \mathbf{R}_{0+} \quad (45)$$

which becomes identical to N if $\vartheta = \varepsilon_0 + \mu + \omega$. □

Note that by taking into account that $\varepsilon_0 \geq \dfrac{\nu N_1}{2N_1 - N_2}$ in Theorem 6, it follows that

$$\left[\dfrac{\nu N_1^2}{(2N_1 - N_2)(\vartheta - \mu - \omega)}, \dfrac{\nu N_1 N_2}{(2N_1 - N_2)(\vartheta - \mu - \omega)}\right] \subset \left[\dfrac{\varepsilon_0 N_1}{\vartheta - \mu - \omega}, \dfrac{\varepsilon_0 N_2}{\vartheta - \mu - \omega}\right]$$

Note also that the total population through time is independent of the vaccination strategy so that it is independent of the ideal vaccination objective constraint $V: \mathbf{R}_{0+} \to \mathbf{R}_+$ as a result. For instance, in a biological war, the objective would be to increase the numbers of the infected plus infectious population for all time. For that purpose, the appropriate vaccination strategy is negative. In the case of the controlled system discussed in Theorem 6. Since the immune population tracks the totals population or a part of it $R^*(t)$, the sum of the remaining populations tracks the difference $N(t) - R^*(t)$ as $t \to \infty$. The case of tracking under constant population; i.e. $\rho = \nu - \mu = 0$ is now discussed. Its proof is direct as a particular case of the proof of Theorem 6 by taking special functions h and g in the vaccination controller.

**Corollary 2**. Assume that $\rho = \nu - \mu = 0$ so that $N(t) = N(0) = N$. Then the following properties hold:

(i) $R(t) = h(t)N$; $\forall t \in \mathbf{R}_{0+}$ if $h(t) = \frac{1}{N}\left( e^{-(\mu+\omega)t} R(0) + \frac{\varepsilon_0 N \left(1 - e^{-(\mu+\omega)t}\right)}{\mu+\omega} \right)$; $\forall t \in \mathbf{R}_{0+}$

If, in addition $\varepsilon_0 = \mu + \omega$ then $h(t) \to 1$ and $R(t) \to N$ as $t \to \infty$.

(ii) $R(t) = h(t)N$; $\forall t \in [0, T)$ and $R(t) = N$; $\forall t \geq T$ for an arbitrary finite time $T \in \mathbf{R}_+$ being non less than a sufficiently large finite design time $T^* = \frac{1}{\mu+\omega} \ln \frac{R(0)}{N} \in \mathbf{R}_+$ provided that $h(t) \in (0,1]$ for $t(\leq T) \in \mathbf{R}_{0+}$, $\varepsilon_0 = \frac{\mu+\omega}{1 - e^{-(\mu+\omega)T}} \left( N - e^{-(\mu+\omega)T} R(0) \right)$ and

$$g(t) = \frac{1}{\varepsilon}\left( 1 - \frac{\mu+\omega}{\varepsilon_0 \left(1 - e^{-(\mu+\omega)t}\right)} \right) \left(1 - e^{-(\mu+\omega)t} R(0)\right); \forall t \geq T$$

**Outline of proof**: Property (i) follows from Theorem 6 by replacing (43) with an identically zero function g. This guarantees that $R(t) = h(t)N$ for all time and also that $h(t) \to 1$ and $R(t) \to N$ as $t \to \infty$ if $\varepsilon_0 = \mu + \omega$. Property (ii) follows in a similar way as the proof of Theorem 6 from (43) by selecting the given time-varying g(t) (which does not dependent now on the dummy time argument $\tau$) in the vaccination controller which guarantees that $R(t) = h(t)N$; $\forall t \in [0, T)$ and $R(t) = N$; $\forall t \geq T \geq T^*$. □

## 6. Modified vaccination law

The vaccination law can be modified by omitting if possible the saturating action (24) so as to improve its effectiveness in collapsing the disease propagation. This may be achieved by generating eventually values exceeding unity at a certain time instant while no partial population is negative at such a time instant. Otherwise, resetting to zero of any potentially negative population is used in the epidemic model. This idea is borrowed from the primary intuitive observation that no partial population can be negative at any time in a real illness propagation situation and the model should reflect that property. The modified vaccination law referred to is as follows for any $t \in \mathbf{R}_{0+}$:

$V(t) = V_a(t)$ if $V_a(t) \geq 0$ and $\min(S(t), E(t), I(t), R(t)) \geq 0$
$V(t) = 1$ if $V_a(t) > 1$ and $\min(S(t), E(t), I(t), R(t)) < 0$
$V(t) = 0$ if $V_a(t) < 0$

where the auxiliary vaccination function $V_a(t)$ is given by (25)-(28) with the subsequent eventual resetting of initial conditions action if some partial population becomes negative at any time instant, that is, if $min(S(t), E(t), I(t), R(t)) < 0$ at time $t \in \mathbf{R}_{0+}$:

$S(t^+) = 0$ if $S(t) < 0$

$E(t^+) = 0$ if $E(t) < 0$

$I(t^+) = 0$ if $I(t) < 0$

$R(t^+) = 0$ if $R(t) < 0$

In the next section, the unsaturated control (249-(28) is compared to the above saturated one in this section through a numerical example.

## 7. Simulation results

This section illustrates through simulation examples the theoretical results stated in the previous sections for the SEIR controlled system. The first example in Section 7.1 is concerned with the saturated vaccination law described by equations (24)-(28) while the second one in Section 7.2 is related to the unsaturated modified vaccination law introduced in Section 6. The SEIR model is described by the following parameters:

$$\frac{1}{\mu} = 255 \text{ days}, \quad \frac{1}{\sigma} = 2.2 \text{ days}, \quad \frac{1}{\omega} = 15 \text{ days}, \quad \rho = 0.1, \quad \frac{1}{\nu} = 150 \text{ days} \quad \beta = 1.66 \text{ days}^{-1}$$

while $\gamma = \sigma$. The initial conditions are given by, S(0) = 400, E(0) = 150, I(0) = 250 and R(0) = 200 individuals so that the total population at initial time is N(0) = 1000 individuals. The function $h(t)$ is selected as recommended in Corollary 2 as $h(t) = \frac{1}{N(t)}\left(e^{-ct}R(0) + N(t)(1 - e^{-ct})\right)$ with $1/c = 5$ days. Notice that $h(t) \to 1$ as $t \to \infty$ and, therefore, $R(t) \to N(t)$. The evolution of the model without vaccination is represented in Figure in order to compare this evolution with the ones associated with different vaccination policies.

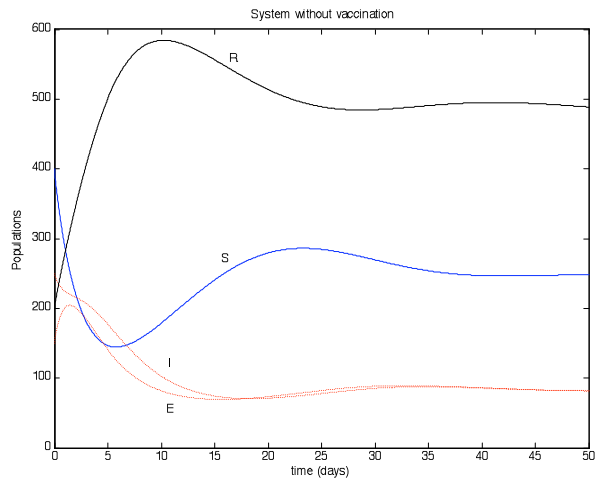

Figure 1. Evolution of populations without vaccination

It is appreciated in Figure 1 that all the populations reach a steady-state value. In particular, there exists a number of infected and infectious individuals, 81 of each, which corresponds in total to an 18% of the total population. The vaccination policies introduced in this work are implemented in order to reduce this percentage of infected and infectious population. Firstly, we consider the saturated vaccination policy given by equations (24)-(28) whose basic feature is that the vaccination effort is restricted to the interval [0,1]. The free controller parameters are selected as $K_r = K_{rd} = 1$. The results are shown in Figures 2 and 3.

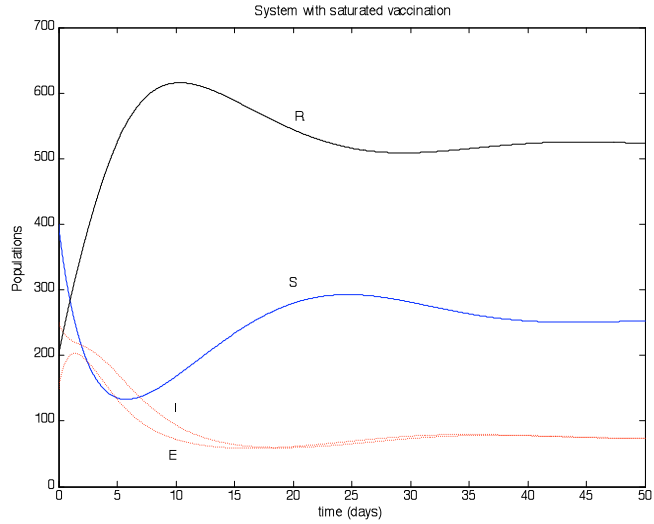

Figure 2. Evolution of the populations with saturated vaccination

Figure 3. Vaccination law for the saturated case

Figures 1 and 2 point out the similarity between the vaccination-free and the saturated vaccination policies: for both of them, the populations reach a steady-state where there exists a non-zero value of infected and infectious people. In particular, each population of infective and infectious is now of 74 individuals which make together a 16% of the total population. This total percentage is slightly smaller to the 18% corresponding to the vaccination free-case. Thus, the reason for the infectious and infective reduction is the application of vaccination with respect to the vaccination-free case. However, the vaccination effort is not large enough to eradicate the illness as it is restricted to the interval [0, 1]. In fact, Now, one considers the saturation to unity of the previous vaccination law is removed under the restriction of all the populations being nonnegative for all time. The results are depicted in Figure 3.

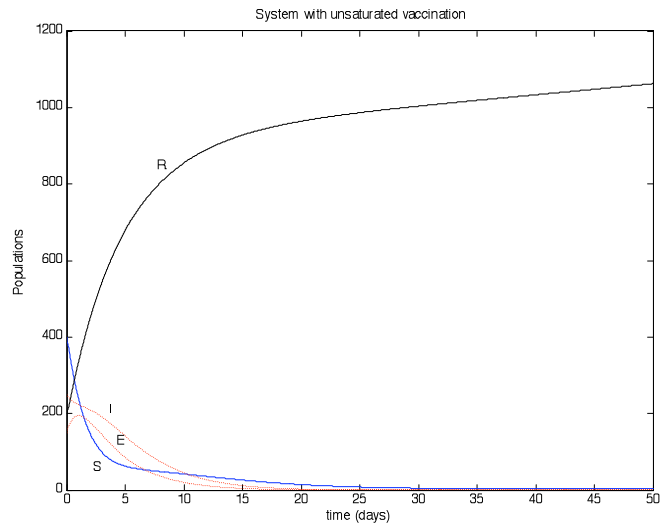

Figure 3. Evolution of the populations under unsaturated vaccination

**Acknowledgements**

The authors thank the Spanish Ministry of Education by its support of this work through Grant DPI2009-07197. He is also grateful to the Basque Government by its support through Grants IT378-10, SAIOTEK SPE07UN04 and SAIOTEK SPE09UN12.

**References**

 [1] M. De la Sen and S. Alonso-Quesada, " A Control theory point of view on Beverton-Holt equation in population dynamics and some of its generalizations", *Applied Mathematics and Computation*, Vol. 199, No. 2, pp. 464-481, 2008.

[2] M. De la Sen and S. Alonso-Quesada, " Control issues for the Beverton-Holt equation in population in ecology by locally monitoring the environment carrying capacity: Non-adaptive and adaptive cases" , *Applied Mathematics and Computation*, Vol. 215, No. 7, pp. 2616-2633, 2009.

[3] M. De la Sen and S. Alonso-Quesada, "Model-matching-based control of the Beverton-Holt equation in Ecology" , *Discrete Dynamics in Nature and Society*, Article number 793512, 2008.

[4] M. De la Sen , "About the properties of a modified generalized Beveron-Holt equation in ecology models", *Discrete Dynamics in Nature and Society*, Article number 592950, 2008.

[5] M. De la Sen, " The generalized  Beverton- Holt equation and the control of populations" , *Applied Mathematical Modelling*, Vol. 32, No. 11, pp. 2312-2328, 2008.

[6] Denis Mollison Editor*, Epidemic Models: Their Structure and Relation to Data*, Publications of the Newton Institute, Cambridge University  Press, 2003.

[7] M. J. Keeling and P. Rohani, *Modeling Infectious Diseases in Humans and Animals*, Princeton University Press, Princeton and Oxford, 2008.

[8] A. Yildirim and Y. Cherruault, "Analytical approximate solution of a SIR epidemic model with constant vaccination strategy by homotopy perturbation method", *Kybernetes*, Vol. 38, No. 9, pp. 1566-1575, 2009.

[9] V.S. Erturk and S. Momani, "Solutions to the problem of prey and predator and the epidemic model via differential transform method", *Kybernetes*, Vol. 37, No. 8 , pp. 1180-1188, 2008.

[10] N. Ortega, L.C. Barros and E. Massad, "Fuzzy gradual rules in epidemiology", *Kybernetes*, Vol. 32, Nos. 3-4, pp. 460-477, 2003.

[11] H. Khan , R.N. Mohapatra, K. Varajvelu and S, J. Liao, " The explicit series solution of SIR and SIS epidemic models", *Applied Mathematics and Computation*, Vol. 215, No. 2 , pp. 653-669, 2009.

[12] X. Y. Song, Y. Jiang and H.M. Wei , " Analysis of a saturation incidence SVEIRS epidemic model with pulse and two time delays", *Applied Mathematics and Computation*, Vol. 214, No. 2 , pp. 381-390, 2009.

[13] T.L. Zhang, J.L. Liu and Z.D. Teng, "Dynamic behaviour for a nonautonomous SIRS epidemic model with distributed delays", *Applied Mathematics and Computation*, Vol. 214, No. 2 , pp. 624-631, 2009.

[14] B. Mukhopadhyay and R. Battacharyya, "Existence of epidemic waves in a disease transmission model with two- habitat population", *International Journal of Systems Science*, Vol. 38, No. 9 , pp. 699-707, 2007.

[15] A. Kelleci and A. Yildirim, "Numerical solution of the system of nonlinear ordinary differential equations arising in kinetic modelling of lactic acid fermentation and epidemic model", *International Journal for Numerical Methods in Biomedical Engineering*", doi: 10.100/cnm.1321.